\nopagenumbers
\input amstex
\voffset-.3in

\font\teneuf=eufm10 \font\seveneuf=eufm7 \font\fiveeuf=eufm5
\font\tenmsy=msbm10 \font\sevenmsy=msbm7 \font\fivemsy=msbm5
\font\tenmsx=msam10 \font\sixmsx=msam6 \font\fivemsx=msam5
\textfont7=\teneuf \scriptfont7=\seveneuf
\scriptscriptfont7=\fiveeuf
\textfont8=\tenmsy \scriptfont8=\sevenmsy
\scriptscriptfont8=\fivemsy
\textfont9=\tenmsx \scriptfont9=\sixmsx
\scriptscriptfont9=\fivemsx

\font\teneuf=eufm10 \font\seveneuf=eufm7 \font\fiveeuf=eufm5
\font\tenmsy=msbm10 \font\sevenmsy=msbm7 \font\fivemsy=msbm5
\font\tenmsx=msam10 \font\sixmsx=msam6 \font\fivemsx=msam5
\textfont7=\teneuf \scriptfont7=\seveneuf
\scriptscriptfont7=\fiveeuf
\textfont8=\tenmsy \scriptfont8=\sevenmsy
\scriptscriptfont8=\fivemsy
\textfont9=\tenmsx \scriptfont9=\sixmsx
\scriptscriptfont9=\fivemsx
\documentstyle{amsppt}
\NoBlackBoxes
\pagewidth{165mm}
\pageheight{260mm}

\ifx\undefined\rom
  \define\rom#1{{\rm #1}}
\fi
\ifx\undefined\curraddr
  \def\curraddr#1\endcurraddr{\address {\it Current address\/}: #1\endaddress}
\fi

\topmatter

\vskip .6in
\title Variants of Schroeder Dissections
\endtitle
\rightheadtext{Schroeder Variants}
\author Leonard M. Smiley\endauthor
\address Len Smiley, Dept. of Mathematical Sciences, UAA, Anchorage AK 99508 \endaddress
\email smiley\@math.uaa.alaska.edu\endemail
\affil Department of Mathematical Sciences, University of Alaska Anchorage, Anchorage, Alaska 99508\endaffil

\abstract Some formulae are given for the enumeration of certain types of dissections of the convex $(n+2)$-gon by non-crossing diagonals. A cataloguing tool, the \it reversive symbol\rm ,\ is proposed.
\endabstract
 
\endtopmatter

 \document
\vskip -10pt
\centerline{July 5, 1999}

\vskip 10pt

\def\l2e{\langle\!\langle}
\def\r2e{\rangle\!\rangle}

\def\Bl2e{\Bigl\langle\!\Bigl\langle}
\def\Br2e{\Bigr\rangle\!\Bigr\rangle}
\noindent
{\bf 1. Dissections and Reversions}
\vskip 0.3 true in
When $m$ diagonals are drawn in a convex $(n+2)$-gon, ($0\leq m\leq n-1$), so that no two diagonals intersect, we call this a \it dissection\rm , and the resulting diagonal-free polygon(s) the \it tiles\rm\  of the dissection. If the vertices of the polygon are labelled, we may ask how many distinct dissections of the $(n+2)$-gon are possible. The answer is an integer sequence known as the solution to Schroeder's Second Problem, or the Schroeder numbers. If in addition we require that $m=n-1$, or, equivalently, that all tiles are triangles, or that there are exactly $n$ tiles, the solution is called the sequence of Catalan numbers. Variant restrictions have been placed on Schroeder dissections and enumerations achieved by a variety of methods (see references in [S1]). A solution may be given by a closed formula or a summation formula for the $n$th term; a recurrence relation defining the sequence; or a closed formula for the generating function of the sequence. An alternative to be emphasized here is the power series whose reversion (or compositional inverse) is the generating function of the solution sequence.
\vskip 0.2 true in
The reasons for this choice are led by two: first, there seems to be a taxonomic advantage to cataloguing restricted dissection counts by the (easily obtained) closed forms of these power series; second, using the formulae of Lagrange and Newton we can often obtain summation formulae for $n$th terms of the solution sequence.

\vskip 0.3 true in
\noindent
{\bf 2. The Triangle-Free Case.} 
\vskip 0.3 true in
Let $a_n$ denote the number of Schroeder dissections of the $(n+2)$-gon having no triangular tiles. If we write $A=A(x)=\sum_{k=0}^\infty a_nx^n$, then the tautological equation $A=1+x^2A^3+x^3A^4+x^4A^5+\cdots$ is easy to see ([S1]). Employing a formal geometric series and a change of variable to $F=xA$, we quickly obtain $x={{F-F^2-F^3}\over{1-F}}:=\alpha (F)$. We call $\alpha$ the \it reversive symbol\rm\  for the sequence $\{ a_n\}$. According to the Lagrange Inversion Formula,
$$a_{n-1}={{1}\over{n}}\Biggl( [t^{n-1}]\Bigl({{t}\over{\alpha (t)}}\Bigr)^n\Biggr),\quad\quad n=1,2,\dots$$
\noindent
and so Newton's formula ('binomial series') gives
$$a_n={{1}\over{n+1}}\sum_{k=0}^{\lceil {{n-1}\over{2}}\rceil}{{n+k}\choose{k}}{{n-k-1}\choose{k-1}}$$
\noindent
where we agree to accept $a_0=1$ in the case of the 2-gon.
\vfill
\eject
\noindent
{\bf 3. Odds and Evens.} 
\vskip 0.2 true in
If we restrict to Schroeder dissections all of whose tiles have an odd number of sides, the reversive symbol is $\alpha (F)={{F-F^2-F^3}\over{1-F^2}}$, and
$$a_n={{1}\over{n+1}}\sum_{k=0}^{\lceil {{n+1}\over{2}}\rceil}{{2n-2k}\choose{n-2k}}{{n-k-1}\choose{k}}$$
follows as before. Here we accept the anomaly that the 2-gon ($n=0$) has $-1$ such dissections! For more detail on this example and the next, see [S2]. 
\vskip 0.2 true in
If all tiles are required to have an even number of sides, the reversive symbol is $\alpha (F)={{F-2F^3}\over{1-F^2}}$. In this case, Newton's formula insists that $a_{2m+1}=0$, as we presumed, and
$$a_{2m}={{1}\over{2m+1}}\sum_{k=0}^m{{2m+k}\choose{k}}{{m-1}\choose{k-1}}.$$
We note here that the sequence $\{ a_{2m}\}$ has been derived in a seemingly different context by Carlitz [C].

\vskip 0.2 true in
\noindent
{\bf 4. Conclusion}
\vskip 0.2 true in
The reversive symbol is useful for the taxonomy of many integer sequences. For the Schroeder numbers themselves we have  $\alpha (F)={{F-2F^2}\over{1-F}}$, and
$$s_n={{1}\over{n+1}}\sum_{k=0}^{n+1}{{2n-k}\choose{n}}{{n-1}\choose{k}}$$
while the simplicity of the Catalan symbol $\alpha (F)=F-F^2$ produces the usual closed formula.
\vskip 0.2 true in
Of course, dissection enumeration is not the only application: the Motzkin numbers, $M_n$, which count the number of ways to place disjoint chords in a circle using $n$ labelled points as endpoints, have symbol $\alpha (F)={{F}\over{1+F+F^2}}={{F-F^2}\over{1-F^3}}$, yielding immediately
$$M_n={{1}\over{n+1}}\sum_{k=0}^{n+1}{{n+1}\choose{k}}{{k}\choose{2k-n-2}}.$$
The author maintains a web-page catalogue of symbols and formulae in the same vein, to which visits and submissions are welcome [S3].
\vskip 0.2 true in
\centerline{Acknowledgements}
\vskip 0.2 true in
I wish to thank D. Zeilberger for the use of EKHAD, and for his E-mailed encouragement to publicize the new dissection counts. Thanks also to Neil Sloane for his interest.   

\enddocument